\documentclass[twoside,10pt]{article}%
\usepackage{amssymb}
\usepackage{amsfonts}
\usepackage{amsmath}
\usepackage{graphicx}%
\providecommand{\U}[1]{\protect\rule{.1in}{.1in}}
 \arraycolsep=1.5pt 

\newtheorem{theorem}{Theorem}

\newtheorem{proposition}[theorem]{Proposition}
\newtheorem{remark}[theorem]{Remark}

\newenvironment{proof}[1][Proof]{\noindent\textbf{#1.} }{\ \rule{0.5em}{0.5em}}
\begin{document}

\title{\textbf{Global Existence of Solutions to the Compressible Euler Equations with Time-dependent Damping and Logarithmic State Equation}}
\author{\textsc{Ka Luen Cheung} and \textsc{Sen Wong\thanks{Corresponding Author, E-mail address: senwongsenwong@yahoo.com.hk}}\\\textit{Department of Mathematics and Information Technology,}\\\textit{The Education University of Hong Kong,}\\\textit{10 Lo Ping Road, Tai Po, New Territories, Hong Kong}}
\date{v55Revised on 28-Jan-2018}
\maketitle
\begin{abstract} In mathematical physics, the pressure function is determined by the equation of state. There are two existing barotropic state equations: the state equation for polytropic gas with adiabatic index greater than or equal to $1$ and the state equation for generalized Chaplygin gas in cosmology. In this paper, a logarithmic pressure is derived naturally with the coexistence of the two existing state equations through an equivalent symmetric hyperbolic transformation. On the study of the logarithmic pressure, global existence of solutions with small initial data to the one-dimensional compressible Euler equations with time-dependent damping is established.
\


\


\end{abstract}
\section{Preliminary}
Consider the following system of equations in a pure mathematical sense.
\begin{align}
    \left\{
\begin{array}
[c]{rl}%
\rho_t +\nabla\cdot(\rho u) & = 0,\\
\rho\lbrack u_{t}+(u\cdot\nabla)u] +\nabla p&= 0,
\end{array}
\right.\label{v45_e1}
\end{align}
where $\rho=\rho(t,x):[0,\infty)\times \mathbb{R}^N\to[0,\infty)$, $u=u(t,x):[0,\infty)\times\mathbb{R}^N\to\mathbb{R}^N$ and $p=p(\rho):[0,\infty)\times\mathbb{R}^N\to\mathbb{R}$ is a function of $\rho$.

In the sense of pure mathematics, compare $(\ref{v45_e1})$ with the following symmetrical hyperbolic system:
\begin{align}
    \left\{
\begin{array}
[c]{rl}%
v_t+B\nabla\cdot u&=-u\cdot \nabla v-\displaystyle\frac{A}{2}v\nabla\cdot u,\\
u_t+B\nabla\cdot v&=-u\cdot\nabla u-\displaystyle\frac{A}{2}v\nabla\cdot v
\end{array}
\right.\label{v45_e2}
\end{align}
for some constants $A\neq0$ and $B\in\mathbb{R}$. Here, $u=u(t,x):[0,\infty)\times\mathbb{R}^N\to\mathbb{R}^N$ and $v=v(t,x):[0,\infty)\times\mathbb{R}^N\to\mathbb{R}$. Then, the following proposition is reached.
\begin{proposition}
With \begin{align}
    v=\frac{2}{A}\left(\sqrt{{p}'(\rho)}-B\right)\label{v51_e3}
\end{align}
and $A\neq0$, system $(\ref{v45_e1})$ is equivalent to system $(\ref{v45_e2})$ if and only if
\begin{align}
    p(\rho)=\left\{
    \begin{matrix}
    \displaystyle\frac{K_1}{A+1}\rho^{A+1}+K,\text{ for }A\neq1,\\
    K_1\ln\rho+K,\text{ for }A=1,
    \end{matrix}
    \right.\label{v45_e4}
\end{align}
where $K_1>0$ and $K$ are arbitrary constants.
\end{proposition}\label{v45_p1}
\begin{remark}
$\sqrt{{p}'(\rho)}$ is well-defined as one can consider complex solutions. 
\end{remark}
\begin{proof}
The proposition can be verified by direct manipulations. We omit the details. 
\end{proof}
\section{Interpretation of Proposition $\ref{v45_p1}$}\label{v48_s2}
The solutions in $(\ref{v45_e4})$ include the following two special cases.
\begin{align}
\text{Case 1: }A>0 \text{ and }A\neq1,\text{ then } p(\rho)=\frac{1}{\gamma}\rho^\gamma,\label{v46_e5}
\end{align}
where $\gamma:=A+1>1$
and
\begin{align}
\text{Case 2: }-2\leq A<-1,\text{ then } p(\rho)=-\frac{1}{\gamma}\rho^{-\gamma},\label{v46_e6}
\end{align}
where $\gamma:=-A-1\in(0,1]$.

Now, note that system $(\ref{v45_e1})$ is in fact the original  $N$-dimensional compressible isentropic Euler equations in fluid mechanics, where $\rho$ and $u$ are the density and the velocity of the fluid respectively. $p=p(\rho)$ is the pressure function determined by a barotropic equation of state. (The term barotropic means that the pressure can be expressed as a function of the density.)

The first equation in $(\ref{v45_e1})$ is derived from the mass conservation law while the second equation in $(\ref{v45_e1})$ is a result of the momentum conservation law. 

For Euler equations in their primitive form, $p$ is given by the gamma law:
\begin{align}
    p=\frac{1}{\gamma}\rho^\gamma, \quad \gamma\geq1.\label{v46_e7}
\end{align}
This ``classical" Euler equations are one of the most fundamental equations in fluid dynamics. Many interesting fluid dynamic phenomena can be described through system $(\ref{v45_e1})$ whose pressure is determined by the gamma law $\cite{v25_r1,v25_r2}$. Moreover, system $(\ref{v45_e1})$ with pressure $(\ref{v46_e7})$ is also the special case of the noted Navier-Stokes equations, whose problem of whether there is a formation of singularity is still open and long-standing.

On the other hand, when the state equation is given by
\begin{align}
    p=-\frac{1}{\gamma}\rho^{-\gamma},\quad 0<\gamma\leq1,\label{v46_e8}
\end{align}
system $(\ref{v45_e1})$ is called the Euler equations for generalized Chaplygin Gas (GCG).

Since 1998, type Ia supernova (SNIa) observations $\cite{v52_r1}$ have shown that our universe has entered into a phase of accelerating expansion. During these years from that time, many additional observational results, including current Cosmic Microwave Background anisotropy measurement $\cite{v52_r2}$, and the data of the Large Scale Structure $\cite{v52_r3}$, also strongly support this suggestion. And these cosmic observations indicate that baryon matter component is about 4$\%$ for total energy density, and about 96$\%$ energy density in universe is invisible. Considering the four-dimensional standard cosmology, this accelerated expansion for universe predicts that dark energy (DE) as an exotic component with negative pressure is filled in universe. And it is shown that DE takes up about two-thirds of the total energy density from cosmic observations. The remaining one third is dark matter (DM). In theory mounting DE models have already been constructed. But there exists another possibility: that the invisible energy component is a unified dark fluid. i.e. a mixture of dark matter and dark energy. The GCG model was introduced in $\cite{v52_r5}$ and developed in $\cite{v52_r4}$ as a unification of dark matter and dark energy.

Now, one can see that the derivations of $(\ref{v46_e5})$ $\&$ $(\ref{v46_e6})$ and $(\ref{v46_e7})$ $\&$ $(\ref{v46_e8})$ are independent. In other words, the state equations $(\ref{v46_e7})$ with $\gamma>1$ and $(\ref{v46_e8})$ which originate from fluid mechanic and cosmology are \textbf{rediscovered} from the equivalent transformation $(\ref{v45_e2})$ without any physical knowledge. However, in addition to $(\ref{v46_e5})$ and $(\ref{v46_e6})$, there is a third possibility in $(\ref{v45_e4})$, namely,
\begin{align}
    p=K_1\ln\rho+K.\label{v46_e9}
\end{align}
We shall call this the Logarithmic Pressure.

Without loss of generality, one may set $K=0$. The following are some highlights of the logarithmic pressure $(\ref{v46_e9})$.\\
1)	It takes both positive and negative values while $(\ref{v46_e7})$ and $(\ref{v46_e8})$ take only positive and negative values respectively.\\
2) As $\ln \rho$ is undefined when $\rho=0$, one has to consider non-vacuum initial data for the system. (The same is true for $(\ref{v46_e8})$).\\
3) $(\ref{v46_e9})$ has the same asymptotic behaviors as $(\ref{v46_e7})$ when $\rho$ tends to infinity and $(\ref{v46_e8})$ when $\rho$ tends to zero.\\
4) $(\ref{v46_e9})$ satisfies the requirement that ${p}'(\rho)>0$ which is the square of sound speed and is required to be positive by a fundamental thermodynamic assumption.

In a broader sense of argument, logarithmic interaction terms like $\log|x-y|$ are common in the interaction energy for planar electrostatics and random matrices. Expressions such as $\rho\log\rho$ are common in information theory and the thermodynamic entropy. In some literature on the semi-geostrophic equations in meteorology \cite{v25_r4}, convexity conditions are imposed on the modified pressure. It would be significant to address which differential equations have long-term solutions and which blow up when the pressure is logarithmic, the marginal case for convexity.

This section is closed by the following two remarks.
\begin{remark}
In $\cite{v47_r1}$, the authors established a shock formation result for the two dimensional (i.e. $N=2$) system $(\ref{v45_e1})$ with state equation $(\ref{v46_e7})$. Within the arguments, a change of variable of $\rho$, namely, $\ln\rho$ was introduced. We remark that $(\ref{v46_e9})$ is not a change of variable but an introduction of possible state equation with potential applications and physical interpretation for the general $N$-dimensional system $(\ref{v45_e1})$. 
\end{remark}
\begin{remark}
In $\cite{v35_r1}$, aiming at proving the global existence of solutions for system $(\ref{v45_e1})$ with pressure $(\ref{v46_e7})$ in the presence of linear damping, the authors considered a special case of the transformation $(\ref{v45_e2})$, where $A=\gamma-1$,  $B=\sqrt{A_1\gamma\bar{\rho}^{\gamma-1}}$ with a positive constant $A_1$ and background density $\bar{\rho}$. In this paper, by investigating the transformation $(\ref{v45_e2})$, which is the primary concern itself, the state equation $(\ref{v46_e9})$ is derived naturally with the coexistence of the two existing state equations.
\end{remark}

\section{Main Result}
On the study of the Logarithmic Pressure, we consider the compressible isentropic Euler equations $(\ref{v45_e1})$ with time-dependent damping which can be expressed as follows.
\begin{equation}
\left\{
\begin{array}
[c]{rl}%
{\normalsize \rho}_{t}{\normalsize +\nabla\cdot(\rho \mathbf{u})} & {\normalsize =}%
{\normalsize 0},\\
\rho\lbrack \mathbf{u}_{t}+(\mathbf{u}\cdot\nabla)\mathbf{u}]{\normalsize +\nabla}p+\displaystyle s(t)\rho\mathbf{u} & {\normalsize =}%
0\text{,}%
\end{array}
\right.  \label{v1_e1}%
\end{equation}
where $s(t)\rho\mathbf{u}$ with $s(t)>0$ is the time-dependent damping term.

As system $(\ref{v1_e1})$ with the logarithmic pressure $(\ref{v46_e9})$ can be transformed to a symmetric hyperbolic system with damping, It follows that finite propagation speed property for that system holds. Furthermore, it is shown in Proposition $\ref{v29_p2}$ that damping can prevent blowup for system $(\ref{v1_e1})$ with the logarithmic pressure $(\ref{v46_e9})$ and small initial data. Hence, it is expected that the logarithmic pressure will share other mathematical and possible physical properties with the two existing pressures.

For system $(\ref{v1_e1})$ with $s(t)=0$ and state equation $(\ref{v46_e7})$, the first general breakdown result for the compressible Euler equations in three spatial dimensions was achieved by Sideris \cite{v41_r2} for adiabatic index $\gamma$ greater than $1$. In \cite{v41_r3}, the author established the finite propagation speed property (F.P.S.P.) by which the followers are able to apply the integration techniques to study the blowup phenomena of system $(\ref{v1_e1})$ when $s(t)=0$. While the global existence result for system $(\ref{v1_e1})$ with $s(t)=0$ is still open, readers may refer to \cite{v35_r6,v35_r7,v35_r8,v35_r9,v41_r3,v41_r4,v41_r5,v41_r6,v41_r7} for the blowup results of system $(\ref{v1_e1})$ without damping.

Subsequently, Sideris et al \cite{v35_r1} showed system $(\ref{v1_e1})$ with constant damping $s(t)=s>0$ and state equation $(\ref{v46_e7})$ still enjoys the finite propagation speed property for $N=3$. Moreover, the authors in \cite{v35_r1} proved the global existence result under the condition that the initial energy functional is small enough.

For system $(\ref{v1_e1})$ with $s(t)=0$ and state equation $(\ref{v46_e8})$, the authors in $\cite{v17_rr1}$ considered the two-dimensional gas expansion problem and proved the global existence of solutions to the expansion problem of a wedge of gas into vacuum with the half angel  $\theta\in(0,\frac{\pi}{2})$ for the generalized Chaplygin gas after obtaining some prior estimates. In \cite{v21_rr1}, the finite propagation speed property for the $N$-dimensional system $(\ref{v1_e1})$ was established. Moreover, by deriving a differential inequality of Sideris' type, it was shown that any $C^1$ solution in a designed non-empty space blows up on finite time provided that the initial functional is sufficiently large.

For $s(t)=\frac{\mu}{(1+t)^\lambda}$, the author in \cite{v41_r1,v35_r3} proved the blowup and global existence results for system $(\ref{v1_e1})$ with state equation $(\ref{v46_e7})$ and $\lambda=1$ in the $1$-dimensional case while the authors in \cite{v35_r2} showed the blowup results for $\lambda\geq1$ and global existence results for $0\leq\lambda\leq1$ for the same system with the same state equation in the $2$ and $3$-dimensional cases.

In this article, we consider system $(\ref{v1_e1})$ with
\begin{align}
    s(t)=\displaystyle\frac{\mu}{1+t}.
\end{align}
The standard approach of energy estimates is applied to obtain the global existence result. Our main contribution is that the result illustrates the mathematical properties of the compressible Euler equations with logarithmic pressure.

For $N=1$, $(\ref{v1_e1})$ becomes
\begin{align}
\left\{
\begin{matrix}
\rho_t+(\rho u)_x=0,\\
(\rho u)_t+(\rho u^2)_x+p_x+\displaystyle\frac{\mu}{1+t}\rho u=0
\end{matrix}\label{v4_e2}
\right.
\end{align}
with the initial data
\begin{equation}\label{v5_e3}
\left\{
\begin{matrix}
(\rho(0,x),u(0,x))=\left(\bar{\rho}+\varepsilon\rho_{0}(x),\varepsilon u_{0}(x)\right),\\
\text{supp}(\rho_0,u_{0})\subseteq\{x:|x|\leq R\}\text{,}%
\end{matrix}
\right.
\end{equation}
for some positive constants $\bar{\rho}$ and $R$. Here $\rho_0$ and $u_0\in H^m(\mathbb{R})$, the Sobolev space with its norm
\begin{align}
    \|f\|_m:=\sum_{k=0}^m\|\partial_x^kf\|_{L^p}.
\end{align}
\begin{remark}
Without loss of generality, one may set $\bar{\rho}=1$.
\end{remark}
We are ready to present the main result, whose proof is given in Section \ref{v42_s3}.
\begin{proposition}\label{v29_p2}
Suppose $(\rho_0,u_0)\in H^m(\mathbb{R})$, $m\geq3$ and $\mu>2$. Then, there exists a unique global classical solution $(\rho(x,t),u(x,t))$ of system $(\ref{v4_e2})$ with state equations $(\ref{v46_e7})$, $(\ref{v46_e8})$ and $(\ref{v46_e9})$. Moreover, the following estimate holds.
\begin{align}
E_m^2(t)+L_m(t)\leq C E_m^2(0),
\end{align}
for some constant $C>0$. Here, $E_m$ and $L_m$ are defined in $(\ref{v5_e16})$ and $(\ref{v40_e22})$ respectively.
\end{proposition}

\section{Proof of Proposition \ref{v29_p2}}\label{v42_s3}
As in $(\ref{v51_e3})$, set
\begin{align}
v:=\frac{2}{K_0}\left(\sqrt{{p}'(\rho)}-\sigma\right),
\end{align}
where $A$ and $B$ in $(\ref{v51_e3})$ are replaced respectively by $K_0$ and
\begin{align}
\sigma:=\sqrt{K_1}.
\end{align}
It is clear that $K_0\neq0,-1$ and $K_0=-1$ correspond to the state equations $(\ref{v46_e7})$ \& $(\ref{v46_e8})$ and $(\ref{v46_e9})$ respectively.

As in $(\ref{v45_e2})$, $(\ref{v4_e2})_1$ and $(\ref{v4_e2})_2$ are transformed to
\begin{align}
\left\{\begin{matrix}
\displaystyle v_t+\sigma u_x&=-uv_x-\frac{K_0}{2}vu_x,\\
\displaystyle u_t+\sigma v_x+\frac{\mu}{1+t}u&=-uu_x-\frac{K_0}{2}vv_x.
\end{matrix}\right.\label{v6_e8}
\end{align}
with initial data
\begin{align}
(v(0,x), u(0,x))=\varepsilon(v_0(x),u_0(x)).
\end{align}
It follows that the system has the finite propagation speed property for the general form of pressure. 

From $(\ref{v6_e8})_1$ and $(\ref{v6_e8})_2$, one has
\begin{align}
v_{tt}-\sigma v_{xx}+\frac{\mu}{1+t}v_t=Q(v,u),\label{v6_e11}
\end{align}
where $Q=Q_1+Q_2+Q_3$ with
\begin{align}
Q_1&=\frac{\mu}{1+t}\left(-uv_x-\frac{K_0}{2}vu_x\right),\label{v6_e12}\\
Q_2&=-\partial_t\left(uv_x+\frac{K_0}{2}vu_x\right),\\
Q_3&=\sigma\partial_x\left(uu_x+\frac{K_0}{2}vv_x\right).\label{v6_e13}
\end{align}
Define
\begin{align}
E_m(T):=\sup_{0<t<T}\{\|(1+t)v_t\|_{m-1}^2+\|(1+t)v_x\|_{m-1}^2+\|(1+t)u_x\|_{m-1}^2+\|v\|^2+\|u\|^2\}^{1/2}\label{v5_e16}
\end{align}
and
\begin{align}
L_m(t):=\int_0^t\left\{(1+\tau)\left(\|v_\tau\|_{m-1}^2+\|v_x\|_{m-1}^2+\|u_x\|_{m-1}^2\right)+\frac{\|u\|^2}{(1+\tau)}\right\}d\tau.\label{v40_e22}
\end{align}


Suppose
\begin{align}
E_m(T)\leq M\varepsilon.\label{v11_e16}
\end{align}
If one can show that
\begin{align}
E_m(T)\leq\frac{1}{2}M\varepsilon,
\end{align}
the the result of global existence follows from the local existence of system $(\ref{v6_e8})$. Here,
$M\geq1$ is independent of $\varepsilon$ and we may assume $\varepsilon$ is small and less than or equal to $1$.

We divide the proof into $7$ steps.

\textbf{Step 1.} Multiply $(\ref{v6_e11})$ by $(1+t)^2v_t$, we have
\begin{align}
\frac{1}{2}\partial_t\left[(1+t)^2v_t^2\right]+(\mu-1)(1+t)v_t^2-\sigma^2(1+t)^2v_tv_{xx}=(1+t)^2v_tQ.\label{v6_e21}
\end{align}
Integrate $(\ref{v6_e21})$ over $[0,t]\times \mathbb{R}$, apply integration by parts and the F.P.S.P., we have
\begin{align}
&\frac{1}{2}\|(1+t)v_t\|^2+\frac{1}{2}\sigma^2\|(1+t)v_x\|^2-\sigma^2\int_0^t(1+\tau)\|v_x\|^2d\tau+(\mu-1)\int_0^t(1+\tau)\|v_\tau\|^2d\tau\label{v6_e22}\\
&=\frac{1}{2}\|v_t(0)\|^2+\frac{1}{2}\sigma^2\|v_x(0)\|^2+\int_0^t\int_\mathbb{R}(1+\tau)^2v_\tau Qdxd\tau.\label{v6_e23}
\end{align}
\textbf{Step 2.} Multiply $(\ref{v6_e11})$ by $(1+t)v$. Then, we have
\begin{align}
\partial_t\left[(1+t)vv_t\right]+\frac{\mu-1}{2}\partial_t(v^2)+\sigma^2(1+t)vv_{xx}-(1+t)(v_t)^2=(1+t)vQ.\label{v6_e24}
\end{align}
Integrate $(\ref{v6_e24})$ over $[0,t]\times \mathbb{R}$, apply integration by parts and the F.P.S.P., we have
\begin{align}
&\int_\mathbb{R}(1+t)vv_tdx+\frac{\mu-1}{2}\|v\|^2+\sigma^2\int_0^t(1+\tau)\|v_x\|^2d\tau-\int_0^t(1+\tau)\|v_\tau\|^2d\tau\label{v6_e25}\\
&=\int_\mathbb{R}(vv_t)(0)dx+\frac{\mu-1}{2}\|v(0)\|^2+\int_0^t\int_\mathbb{R}(1+\tau)vQdxd\tau.\label{v6_e26}
\end{align}
The first two terms of $(\ref{v6_e23})$ and the second term of $(\ref{v6_e26})$ are easy to handle as their sum is less than or equal to $CE^2_1(0)$, for some positive $C$. By the AM-GM inequity, the first term of $(\ref{v6_e26})$ is less than or equal to $CE^2_1(0)$, for some positive $C$. In what follows, $C$, which may change from line to line, will denote a positive constant depending on $\mu$, $K_0$, $K_1$ and $\sigma$ only.

On the other hand, the first two terms of $(\ref{v6_e22})$ and the second term of $(\ref{v6_e25})$ should be kept as they appear in $E^2_1(T)$. Also, the fourth term of $(\ref{v6_e22})$ and the fourth term of $(\ref{v6_e25})$ should be kept as they apper in $L_1(t)$.

Now, we want to keep the term $\displaystyle\int_0^t(1+\tau)\|v_x\|^2d\tau$ as it appears in $L_1$. Thus, it remains to cope with the first term of $(\ref{v6_e25})$.

By AM-GM inequality,
\begin{align}
v(1+t)v_t\geq-\frac{b}{4}v^2-\frac{1}{b}(1+t)^2(v_t)^2
\end{align}
for any positive $b$. Thus,
\begin{align}
\text{the first term of $(\ref{v6_e25})$}\geq-\frac{b}{4}\|v\|^2-\frac{1}{b}\|(1+t)v_t\|^2.\label{v6_e28}
\end{align}
Thus, the first term of $(\ref{v6_e28})$ can be absorbed by the second term of $(\ref{v6_e25})$ if $b$ is small enough; and the second term of $(\ref{v6_e28})$ can be absorbed by the the first term of $(\ref{v6_e22})$. However, we wish to keep the term $\int_0^t(1+\tau)\|v_x\|^2d\tau$. Thus, assuming $\mu>2$, we multiply $(\ref{v6_e22})$ by a constant $a$ and add the result to $(\ref{v6_e25})$. In this way, one should carefully choose $a$ and $b$ to keep all the coefficients positive. After examination,
\begin{align}
a=\frac{\mu}{2(\mu-1)}\text{ and }b=\frac{8(\mu-1)}{\mu-2}
\end{align}
will work. Thus, after $a\times(\ref{v6_e22})$+$(\ref{v6_e25})$ and using $(\ref{v6_e28})$, one obtains
\begin{align}
&\|(1+t)v_t\|^2+\|(1+t)v_x\|^2+\|v\|^2+\int_0^t(1+\tau)\|v_x\|^2d\tau+\int_0^t(1+\tau)\|v_\tau\|^2d\tau\label{v7_e28}\\
&\leq C E^2_1(0)+C\left[\int_0^t\int_\mathbb{R}a(1+\tau)^2v_\tau Qdxd\tau+\int_0^t\int_\mathbb{R}(1+\tau)vQdxd\tau\right], 
\end{align}
where $C$ is a positive constant independent of $M$ and $\varepsilon$.

\textbf{Step 3.} Multiply $(\ref{v6_e8})_2$ by u, take integration over $[0,t]\times \mathbb{R}$, use integration by parts, one has
\begin{align}
&\frac{1}{2}\|u\|^2+\mu\int_0^t\frac{\|u\|^2}{1+\tau}d\tau-\sigma\int_0^t\int_\mathbb{R}vu_xdxd\tau\label{v6_e32}\\
&=\frac{1}{2}\|u(0)\|^2-\int_0^t\int_\mathbb{R}(u^2u_x+\frac{K_0}{2}uvv_x)dxd\tau.
\end{align}
To handle the third term of $(\ref{v6_e32})$, we use $(\ref{v6_e8})_2$ and integration by parts to obtain
\begin{align}
\text{the third term of $(\ref{v6_e32})$ }=\frac{1}{2}\|v\|^2-\frac{1}{2}\|v(0)\|^2+\int_0^t\int_\mathbb{R}uvv_xdxd\tau.
\end{align}
Thus, one has
\begin{align}
&\|u\|^2+\|v\|^2+\int_0^t\frac{\|u\|^2}{1+\tau}d\tau\label{v7_e35}\\
&\leq C \left[E^2_1(0)+\int_0^t\int_\mathbb{R}\left[-(1+\frac{K_0}{2})uvv_x-u^2u_x\right]dxd\tau\right],
\end{align}
for some positive constant $C$.

\textbf{Step 4.} Differentiate $(\ref{v6_e8})_2$ with respect to $x$, multiply it by $(1+t)^2u_x$, integrate it over $[0,t]\times\mathbb{R}$ and apply integration by parts. One has
\begin{align}
&\frac{1}{2}\|(1+t)u_x\|^2+(\mu-1)\int_0^t(1+\tau)\|u_x\|^2d\tau-\sigma\int_0^t\int_\mathbb{R}(1+\tau)^2u_{xx}dxd\tau\label{v7_e37}\\
&=\frac{1}{2}\|u_x(0)\|^2+\int_0^t\int_\mathbb{R}(1+\tau)^2u_x\partial_x(-uu_x-\frac{K_0}{2}vv_x)dxd\tau.\label{v7_e38}
\end{align}
By $(\ref{v6_e8})_1$, 
\begin{align}
\text{the third term of $(\ref{v7_e37})$}&=\int_0^t\int_\mathbb{R}(1+\tau)^2v_x\left[v_{x\tau}+\partial_x(uv_x+\frac{K_0}{2}vu_x)\right]dxd\tau\\
&=\frac{1}{2}\|(1+t)v_x\|^2-\frac{1}{2}\|v_x(0)\|^2-\int_0^t(1+\tau)\|v_x\|^2d\tau+\int_0^t\int_\mathbb{R}(1+\tau)^2v_x\partial_x(uv_x+\frac{K_0}{2}vu_x)dxd\tau.
\end{align}
Hence, $(\ref{v7_e37})$ and $(\ref{v7_e38})$ become
\begin{align}
&\frac{1}{2}\|(1+t)u_x\|^2+\frac{1}{2}\|(1+t)v_x\|^2+(\mu-1)\int_0^t(1+\tau)\|u_x\|^2d\tau-\int_0^t(1+\tau)\|v_x\|^2d\tau\label{v7_e41}\\
&\leq CE^2_1(0)+\int_0^t\int_\mathbb{R}(1+\tau)^2u_x\partial_x(-uu_x-\frac{K_0}{2}vv_x)dxd\tau+\int_0^t\int_\mathbb{R}(1+\tau)^2v_x\partial_x(-uv_x-\frac{K_0}{2}vu_x)dxd\tau.
\end{align}
Note that the negative term (the fourth term) of $(\ref{v7_e41})$ can be absorbed by the forth term of $(\ref{v7_e28})$. Now we add the following.
\begin{align}
2\times\text{(\ref{v7_e28}) in Step 2 }+\text{(\ref{v7_e35}) in Step 3 }+\text{(\ref{v7_e41}) in Step 4}
\end{align}
to have
\begin{align}
E^2_1(T)+L_1(t)&\leq CE^2_1(0)+C\left[|I_1|+|I_2|+|I_3|+|J_1|+|J_2|+|W_1|+|W_2|+|W_3|+|W_4|\right],\label{v17_e41} 
\end{align}
where
\begin{align}
I_1&=\int_0^t\int_\mathbb{R}(1+\tau)v_\tau\left(vv_x+\frac{K_0}{2}vu_x\right)dxd\tau,\\
I_2&=\int_0^t\int_\mathbb{R}v\left(uv_x+\frac{K_0}{2}vu_x\right)dxd\tau,\\
I_3&=\int_0^t\int_\mathbb{R}\left[\left(1+\frac{K_0}{2}\right)uvv_x+u^2u_x\right]dxd\tau,\\
J_1&=\int_0^t\int_\mathbb{R}(1+\tau)v\partial_\tau\left(uv_x+\frac{K_0}{2}vu_x\right)dxd\tau,\\
J_2&=\int_0^t\int_\mathbb{R}(1+\tau)v\partial_x\left(uu_x+\frac{K_0}{2}vv_x\right)dxd\tau,\\
W_1&=\int_0^t\int_\mathbb{R}(1+\tau)^2v_\tau\partial_\tau\left(uv_x+\frac{K_0}{2}vu_x\right)dxd\tau,\\
W_2&=\int_0^t\int_\mathbb{R}(1+\tau)^2v_\tau\partial_x\left(uu_x+\frac{K_0}{2}vv_x\right)dxd\tau,\\
W_3&=\int_0^t\int_\mathbb{R}(1+\tau)^2u_x\partial_x\left(uu_x+\frac{K_0}{2}vv_x\right)dxd\tau,\\
W_4&=\int_0^t\int_\mathbb{R}(1+\tau)^2v_x\partial_x\left(uv_x+\frac{K_0}{2}vu_x\right)dxd\tau.
\end{align}
\textbf{Step 5.} In this step, we will handle $I_k$ and $J_k$.

First, by Sobolev's inequality and $(\ref{v11_e16})$, we have
\begin{align}
\begin{matrix}
\displaystyle\sup_{t,x}\left\{|u|,|v|,|u_x|,|v_x|,|v_t|, (1+t)|v_t|, (1+t)|v_x|, (1+t)|u_x|\right\}\leq CM\varepsilon.
\end{matrix}
\label{v11_e51} 
\end{align}
For $I_1$, note that by $(\ref{v11_e51})$ and AM-GM inequality, one has
\begin{align}
|v_\tau vv_x|+|\frac{K_0}{2}v_\tau vu_x|&\leq CM\varepsilon\left(|v_\tau v_x|+|v_\tau u_x|\right)\\
&\leq CM\varepsilon\left(|v_\tau|^2+|v_x|^2+|u_x|^2\right).
\end{align}
Hence,
\begin{align}
|I_1|\leq CM\varepsilon L_1(t).
\end{align}
For $I_2$, by integration by parts, one has
\begin{align}
I_2=\int_0^t\int_\mathbb{R}(1-K_0)uvv_xdxd\tau.
\end{align}
As
\begin{align}
|uvv_x|&\leq CM\varepsilon|uv_x|\\
&=CM\varepsilon\left|\frac{u}{\sqrt{1+\tau}}\right|\left|\sqrt{1+\tau}v_x\right|\\
&\leq CM\varepsilon\left(\frac{u^2}{1+\tau}+(1+\tau)(v_x)^2\right),
\end{align}
one has
\begin{align}
|I_2|\leq CM\varepsilon L_1(t).
\end{align}
Similarly, one also has
\begin{align}
|I_3|\leq CM\varepsilon L_1(t).
\end{align}
For $J_1$, by integration by parts,
\begin{align}
J_1=\int_0^t\int_\mathbb{R}(1+\tau)\left[-v_\tau v_x u+\left(-1+\frac{K_0}{2}\right)vv_\tau u_x+\left(1-K_0\right)v v_x u_\tau\right]dxd\tau.
\end{align}
Hence,
\begin{align}
|J_1|\leq CM\varepsilon L_1(t)+CM\varepsilon\int_0^t\int_\mathbb{R}(1+\tau)|u_\tau|^2dxd\tau.
\end{align}
Note that by $(\ref{v6_e8})_2$ and $(\ref{v11_e51})$,
\begin{align}
|u_\tau|\leq CM\left(|v_x|+|u_x|+\frac{|u|}{1+\tau}\right).
\end{align}
By AM-GM inequality,
\begin{align}
|u_\tau|^2\leq CM^2\left[|v_x|^2+|u_x|^2+\frac{|u|^2}{(1+\tau)^2}\right].
\end{align}
Hence,
\begin{align}
|J_1|&\leq CM\varepsilon L_1(t)+CM^3\varepsilon L_1(t)\\
&\leq CM^3\varepsilon L_1(t).
\end{align}
For $J_2$, by integration by parts,
\begin{align}
J_2=\int_0^t\int_\mathbb{R}(1+\tau)\left[-uu_xv_x-\frac{K_0}{2}v(v_x)^2\right]dxd\tau.
\end{align}
Thus, one has
\begin{align}
|J_2|\leq CM^3\varepsilon L_1(t).
\end{align}
\textbf{Step 6.} In this step, we shall handle $W_k$.

First, by integration by parts, one has
\begin{align}
W_1&=\int_0^t\int_\mathbb{R}(1+\tau)^2\left[\left(\frac{K_0}{2}-\frac{1}{2}\right)u_x(v_\tau)^2+u_\tau v_xv_\tau\right]dxd\tau\label{v13_e66}\\
&+\frac{K_0}{2}\int_0^t\int_\mathbb{R}(1+\tau)^2vv_\tau u_{x\tau}dxd\tau\label{v13_e67},\\
W_2&=\int_0^t\int_\mathbb{R}(1+\tau)^2\left[v_\tau(u_x)^2-\frac{K_0}{2}vv_xv_{x\tau}\right]dxd\tau\label{v13_e68}\\
&+\int_0^t\int_\mathbb{R}(1+\tau)^2uv_\tau u_{xx}dxd\tau\label{v13_e69},\\
W_3&=\int_0^t\int_\mathbb{R}(1+\tau)^2\left[\frac{1}{2}u_x(u_x)^2\right]dxd\tau\label{v13_e70}\\
&-\frac{K_0}{2}\int_0^t\int_\mathbb{R}(1+\tau)^2v v_x u_{xx}dxd\tau\label{v13_e71},\\
W_4&=\int_0^t\int_\mathbb{R}(1+\tau)^2\left[\left(\frac{K_0}{2}+\frac{1}{2}\right)u_x(v_x)^2\right]dxd\tau\label{v13_e72}\\
&+\frac{K_0}{2}\int_0^t\int_\mathbb{R}(1+\tau)^2vv_xu_{xx}dxd\tau.\label{v13_e73}
\end{align}
Second, the right hand sides of $(\ref{v13_e66})$, $(\ref{v13_e68})$, $(\ref{v13_e70})$ and $(\ref{v13_e72})$ are easy to handle as they only contain at most first order derivatives. The exception is the term that contains $v_{x\tau}$ in $(\ref{v13_e68})$, which will be handled as follows.

Note that
\begin{align}
vv_xv_{x\tau}=\frac{1}{2}v\partial_\tau\left[\left(v_x\right)^2\right].\label{v13_e74}
\end{align}
Hence, the second order derivative will be reduced to first order derivative when integration by parts is applied. More precisely,
\begin{align}
-\frac{K_0}{2}\int_0^t\int_\mathbb{R}(1+\tau)^2vv_xv_{x\tau}dxd\tau&=-\frac{K_0}{2}\int_\mathbb{R}v(1+t)^2(v_x)^2dx+\frac{K_0}{2}\int_\mathbb{R}v(0)(v_x(0))^2dx\\
&+\frac{K_0}{2}\int_0^t\int_\mathbb{R}(1+\tau)(v_x)^2\left[(1+\tau)v_\tau+2v\right]dxd\tau\\
&\leq CM\varepsilon\|(1+t)v_x\|^2+CM\varepsilon E_1^2(0)+CM\varepsilon L_1(t).
\end{align}
Hence,
\begin{align}
\text{ Right hand side of }(\ref{v13_e68})\leq CM\varepsilon\|(1+t)v_x\|^2+CM\varepsilon E_1^2(0)+CM\varepsilon L_1(t).
\end{align}
Similarly, we have
\begin{align}
\text{ Right hand side of }(\ref{v13_e66})&\leq CM\varepsilon L_1(t),\\
\text{ Right hand side of }(\ref{v13_e70})&\leq CM\varepsilon L_1(t),\\
\text{ Right hand side of }(\ref{v13_e72})&\leq CM\varepsilon L_1(t).
\end{align}
It remains to handle $(\ref{v13_e67})$, $(\ref{v13_e69})$, $(\ref{v13_e71})$ and $(\ref{v13_e73})$. The method is similar to that one in $(\ref{v13_e74})$. To be more specific, whenever $v_xv_{xx}$. $v_\tau v_{\tau\tau}$, $v_\tau v_{\tau x}$, $v_x v_{x\tau}$, $u_xu_{xx}$, $u_\tau u_{\tau\tau}$, $u_\tau u_{\tau x}$ or $u_x u_{x\tau}$ occur, one may apply integration by parts to reduce them to first order derivatives. Note that by $(\ref{v6_e8})_1$, one has
\begin{align}
u_{x\tau}&=\frac{-u v_{x\tau}-v_{\tau\tau}-u_\tau v_x}{\sigma+\frac{K_0}{2}v}+\frac{K_0}{2}\frac{(v_\tau)^2+uv_xv_\tau}{\left(\sigma+\frac{K_0}{2}v\right)^2}
\end{align}
and
\begin{align}
u_{xx}&=\frac{-uv_{xx}-v_{x\tau}-u_xv_x}{\sigma+\frac{K_0}{2}v}+\frac{K_0}{2}\frac{u(v_x)^2+v_xv_\tau}{\left(\sigma+\frac{K_0}{2}v\right)^2}.
\end{align}
To handle the denominator $(\sigma+\frac{K_0}{2}v)$, note that
\begin{align}
\sigma+\frac{K_0}{2}v\geq \sigma(1-C_1M\varepsilon)>0
\end{align}
for some positive constant $C_1$.
Hence,
\begin{align}
(\ref{v13_e67})&\leq\frac{CM\varepsilon}{1-C_1M\varepsilon}\|(1+t)v_t\|^2+\frac{CM\varepsilon}{1-C_1M\varepsilon}E_1^2(0)+\frac{CM^3\varepsilon}{(1-C_1M\varepsilon)^2}L_1(t)\\
(\ref{v13_e69})&\leq\frac{CM^2\varepsilon^2}{1-C_1M\varepsilon}\|(1+t)v_x\|^2+\frac{CM^2\varepsilon^2}{1-C_1M\varepsilon}E_1^2(0)+\frac{CM^3\varepsilon}{(1-C_1M\varepsilon)^2}L_1(t)\\
(\ref{v13_e71})&\leq\frac{CM\varepsilon}{1-C_1M\varepsilon}\|(1+t)v_x\|^2+\frac{CM\varepsilon}{1-C_1M\varepsilon}E_1^2(0)+\frac{CM^3\varepsilon}{(1-C_1M\varepsilon)^2}L_1(t)\\
(\ref{v13_e73})&\leq\frac{CM\varepsilon}{1-C_1M\varepsilon}\|(1+t)v_x\|^2+\frac{CM\varepsilon}{1-C_1M\varepsilon}E_1^2(0)+\frac{CM^3\varepsilon}{(1-C_1M\varepsilon)^2}L_1(t).
\end{align}
\textbf{Final Step.} Now, $(\ref{v17_e41})$ becomes
\begin{align}
E_1^2(t)+L_1(t)\leq \frac{C(1-C_1M\varepsilon+M^2\varepsilon)}{1-C_1M\varepsilon}E_1^2(0)+\frac{CM^3\varepsilon}{(1-C_1M\varepsilon)^2}E^2_1(t)+\frac{CM^3\varepsilon}{(1-C_1M\varepsilon)^2}L_1(t).
\end{align}
In general, one has
\begin{align}
E_m^2(t)+L_m(t)\leq \frac{C(1-C_1M\varepsilon+M^2\varepsilon)}{1-C_1M\varepsilon}E_m^2(0)+\frac{CM^3\varepsilon}{(1-C_1M\varepsilon)^2}E^2_m(t)+\frac{CM^3\varepsilon}{(1-C_1M\varepsilon)^2}L_m(t).
\end{align}
Hence, 
\begin{align}
E_m^2(t)+L_m(t)\leq C_2 E_m^2(0),
\end{align}
where
\begin{align}
C_2=\frac{C(1-C_1M\varepsilon)(1-C_M\varepsilon+M^2\varepsilon)}{(1-C_1M\varepsilon)^2-CM^3\varepsilon}.
\end{align}
There exists an $M_0\geq1$ such that for any given $M\geq M_0$, one may choose an $\varepsilon_0<\frac{1}{C_1M}$ such that $C_2\leq C_1\varepsilon_0+1$. Given $E_m^2(0)\leq C_3\varepsilon^2$, set $\frac{M^2}{4}:=\max\{\frac{M_0^2}{4}, (C_1\varepsilon+1)C_3\}$, then
\begin{align}
E^2_m(t)\leq \frac{1}{4}M^2\varepsilon^2.
\end{align}
The proof is completed.





\section{Acknowledgement}
The research was partially supported by the EDUHK Dean's Research Fund 17/18 and small scale grant 18/19.

\end{document}